\documentclass[12pt]{amsart}

\usepackage{amscd}
\usepackage{srcltx}
\usepackage{amsmath}
\usepackage{amsfonts}
\usepackage{verbatim}
\usepackage{amsthm}
\usepackage{amssymb, mathrsfs}
\usepackage{amscd}
\usepackage{latexsym}
\usepackage{array}
\usepackage{enumerate}
\usepackage{longtable}
\usepackage[all]{xypic}
\usepackage[latin1]{inputenc}

\ifx\pdfpagewidth\undefined 
\usepackage[T1]{fontenc}
\else 
\usepackage[OT1]{fontenc} 
\fi 
\usepackage{amsfonts,amssymb,latexsym,longtable,amsmath}

\def\sbs{\subset}
\def\obr{^{-1}}
\def\stm{\setminus}

\def\norm#1{\left\Vert#1\right\Vert}

\def\g{\gamma}
\def\d{\delta}
\def\o{\omega}

\newtheorem{thm}{Theorem}[section]
\newcommand{\bthm}{\begin{thm}} \newcommand{\ethm}{\end{thm}}
\newtheorem{prop}[thm]{Proposition}
\newcommand{\bprp}{\begin{prop}} \newcommand{\eprp}{\end{prop}}
\newtheorem{fact}[thm]{Fact}
\newcommand{\bfct}{\begin{fact}} \newcommand{\efct}{\end{fact}}
\newtheorem{prob}[thm]{Problem}
\newcommand{\bprb}{\begin{prob}} \newcommand{\eprb}{\end{prob}}
\newtheorem{question}[thm]{Question}
\newcommand{\bque}{\begin{question}} \newcommand{\eque}{\end{question}}
\newtheorem{lem}[thm]{Lemma}
\newcommand{\blem}{\begin{lem}} \newcommand{\elem}{\end{lem}}
\newtheorem{claim}[thm]{Claim}
\newcommand{\bclm}{\begin{claim}} \newcommand{\eclm}{\end{claim}}
\newtheorem{cor}[thm]{Corollary}
\newcommand{\bcor}{\begin{cor}} \newcommand{\ecor}{\end{cor}}
\newtheorem{conj}[thm]{Conjecture}
\newcommand{\bcnj}{\begin{conj}} \newcommand{\ecnj}{\end{conj}}
\theoremstyle{definition}
\newtheorem{defn}[thm]{Definition}
\newcommand{\bdfn}{\begin{defn}} \newcommand{\edfn}{\end{defn}}
\newtheorem{spec}[thm]{Specializing}
\newcommand{\bspc}{\begin{spec}} \newcommand{\espc}{\end{spec}}
\theoremstyle{remark}
\newtheorem{rem}[thm]{Remark}
\newcommand{\brem}{\begin{rem}} \newcommand{\erem}{\end{rem}}
\newtheorem{cnv}[thm]{Convention}
\newcommand{\bcnv}{\begin{cnv}} \newcommand{\ecnv}{\end{cnv}}
\newtheorem{exam}[thm]{Example}
\newcommand{\bexm}{\begin{exam}} \newcommand{\eexm}{\end{exam}}
\newcommand{\bpf}{\begin{proof}} \newcommand{\epf}{\end{proof}}

\newcommand{\C}{\mathbb C}

\newcommand{\sK} {{\mathcal K}}

\newcommand{\sN} {{\mathcal N}}

\renewcommand{\phi}{\varphi}
\renewcommand{\theta}{\vartheta}

\newcommand{\gep}{{\epsilon}}

\newcommand{\s}{{\sigma}}

\begin{document}
\title{The dual space of precompact groups}

\author[M. Ferrer]{M. Ferrer}
\address{Universitat Jaume I, Instituto de Matem\'aticas de Castell\'on,
Campus de Riu Sec, 12071 Castell\'{o}n, Spain.}
\email{mferrer@mat.uji.es}

\author[S. Hern\'andez]{S. Hern\'andez}
\address{Universitat Jaume I, INIT and Departamento de Matem\'{a}ticas,
Campus de Riu Sec, 12071 Castell\'{o}n, Spain.}
\email{hernande@mat.uji.es}

\author[V. Uspenskij]{V. Uspenskij}
\address{Department of mathematics, 321 Morton Hall, Ohio University, Athens, Ohio 45701, USA}
\email{uspenski@ohio.edu}

\thanks{ The first and second listed
authors acknowledge partial support by the Spanish Ministry of
Science, grant MTM2008-04599/MTM}

\begin{abstract}
For any topological group $G$ the dual object $\widehat G$ is defined as the set of equivalence classes
of irreducible unitary representations of $G$ equipped with the Fell topology. If $G$ is compact, $\widehat G$ is
discrete. In an earlier paper we proved that $\widehat G$ is discrete for every metrizable precompact group, 
i.e.~a dense subgroup of a compact metrizable group. We generalize this result to the case when $G$
is an almost metrizable precompact group.
\end{abstract}

\thanks{{\em 2010 Mathematics Subject Classification.} Primary 43A40. Secondary 22A25, 22C05, 22D35, 43A35, 43A65, 54H11\\
{\em Key Words and Phrases:} compact group, precompact group, representation,
Pontryagin--van Kampen duality, compact-open topology, Fell dual space, Fell topology, 
Kazhdan property (T)}


\date{22 December 2011}

\maketitle \setlength{\baselineskip}{24pt}
\setlength{\parindent}{1cm}


\section{Introduction}
\label{s:intro}

For a topological group $G$ let $\widehat G$ be the set of equivalence classes
of irreducible unitary representations of $G$. The set $\widehat G$ can be equipped
with a natural topology, the so-called Fell topology (see Section~\ref{s:prelim} for a definition).

A topological group $G$ is {\em precompact}
if it is isomorphic (as a topological group) to a subgroup of a compact
group $H$ (we may assume that $G$ is dense in $H$).
If $H$ is compact, then $\widehat H$ is discrete. If $G$ is a dense
subgroup of $H$, the natural mapping $\widehat H\to \widehat G$ is a bijection but
in general need not be a homeomorphism. Moreover, for every countable
non-metrizable precompact group $G$ the space $\widehat G$ is not discrete \cite[Theorem 5.1]{FHU}, 
and every non-metrizable compact group $H$ has a dense subgroup $G$ such that $\widehat G$
is not discrete \cite[Theorem 5.2]{FHU}. 
(The Abelian case was considered in \cite{comractri:04,dik_sha:jmaa_det,hermactri}).
On the other hand, if $G$ is a precompact metrizable group, then $\widehat G$ is discrete
\cite[Theorem 4.1]{FHU}. (The Abelian case was considered in \cite{aus,chasco}).
The aim of the present paper is to generalize this result to the almost metrizable case:
$\widehat G$ is discrete for every almost metrizable
precompact topological group $G$. A topological group $G$ is {\em almost metrizable} if it has a compact
subgroup $K$ such that the quotient space $G/K$ is metrizable. According to Pasynkov's theorem
\cite[4.3.20]{ArhTk}, a topological group is almost metrizable if and only if it is feathered in the sense
of Arhangel'skii.

%
We reduce the almost metrizable case to the metrizable case considered in \cite[Theorem 4.1]{FHU}.

\section{Preliminaries: Fell topologies}
\label{s:prelim}

All topological spaces and groups that we consider are assumed to be Hausdorff.
For a (complex) Hilbert space $\mathcal{H}$ the unitary group $U(\mathcal{H})$
of all linear isometries of $\mathcal{H}$ is equipped with the strong operator topology
(this is the topology of pointwise convergence). With this topology, $U(\mathcal{H})$ is a topological group.

A {\em unitary representation} $\rho$ of the to\-po\-lo\-gi\-cal group $G$
is a continuous homomorphism $G\to U(\mathcal{H})$, where $\mathcal{H}$ is a complex Hilbert space.
A closed linear subspace $E\subseteq \mathcal H$ is an \emph{invariant} subspace
for $\mathcal S\subseteq U(\mathcal{H})$ if $ME\subseteq E$ for all
$M\in \mathcal S$.
If there is a closed subspace $E$ with $\{0\}\subsetneq E\subsetneq \mathcal H$
which is invariant for $\mathcal S$, then $\mathcal S$ is called
\emph{reducible}; otherwise $\mathcal S$ is \emph{irreducible}.
An \emph{irreducible representation} of $G$ is a 
unitary representation $\rho$ such that 
$\rho(G)$ is irreducible.

If $\mathcal{H}=\C^n$, we identify $U(\mathcal{H})$ with the {\em unitary group of order $n$},
that is, the compact Lie group of all complex $n\times n$ matrices $M$ for which $M^{-1}=M^{*}$.
We denote this group by $\mathbb{U}(n)$.

Two unitary representations $\rho:G\to U(\mathcal{H}_1)$ and $\psi: G\to U(\mathcal{H}_2)$
are {\it equivalent} 
if there exists a Hilbert space isomorphism $M:\mathcal {H}_1\to \mathcal H_2$
such that $\rho(x)=M^{-1}\psi(x)M$ for all
$x\in G$. The {\em dual object} of a topological group $G$ is the set $\widehat G$ of
equivalence classes of irreducible unitary representations of $G$.

If $G$ is a precompact group, the Peter-Weyl Theorem (see \cite{hof_mor:compact_groups}) implies that
all irreducible unitary representation of $G$ are finite-dimensional and determine an embedding of $G$ into the product
of unitary groups $\mathbb{U}(n)$.

If $\rho:G\to U(\mathcal{H})$ is a unitary representation, a complex-valued function $f$ on $G$ is called
a {\em function of positive type} (or {\em positive-definite function}) {\em associated with} $\rho$ if there exists a vector
$v\in \mathcal{H}$ such that $f(g)=(\rho(g)v, v)$
(here $(\cdot ,\cdot)$ denotes the inner product in $\mathcal{H}$). We denote by
$P_\rho'$ the set of all functions of positive type associated with $\rho$.
Let $P_\rho$ be the convex cone generated by $P_\rho'$, that is, the set of sums of elements of $P_\rho'$.

Let $G$ be a topological group, $\mathcal R$ a set of equivalence classes of unitary representations of $G$.
The \emph{Fell topology} on $\mathcal R$ is defined as follows: a typical neighborhood of $[\rho]\in \mathcal R$ has the form
$$
W(f_1, \cdots, f_n, C, \gep)=\{[\s]\in\mathcal R :
\, \exists g_1, \cdots, g_n\in P_\s
\  \forall x\in C  \
|f_i(x)-g_i(x)| <\gep\},
$$
where $f_1, \cdots, f_n\in P_\rho$ (or $\in P_\rho'$), $C$ is a compact subspace of $G$, and $\gep>0$.
In particular, the Fell topology is defined on the dual object $\widehat G$. If $G$ is locally compact,
the Fell topology on $\widehat G$ can be derived from the Jacobson topology on the primitive ideal space of $C^*(G)$,
 the $C^*$-algebra of $G$
\cite[section 18]{dixmier}, \cite[Remark F.4.5]{bhv:book}.

Every onto homomorphism $f:G\to H$ of topological groups gives rise to a continuous injective dual map $\hat f:\widehat H\to \widehat G$.
A mapping $h:X\to Y$ between topological spaces is {\em compact-covering}
if for every compact set $L\sbs Y$ there exists a compact set $K\sbs X$ such that $h(K)=L$. 

\blem
\label{l:ccov}
If $f:G\to H$ is a compact-covering onto homomorphism of topological groups, the dual map $\hat f:\widehat H\to \widehat G$
is a homeomorphic embedding.
\elem

\bpf
This easily follows from the definition of Fell topology.
\epf

Let $\pi$ be a unitary representation of a topological group $G$ on a Hilbert space $\mathcal H$.
Let $F\subseteq G$ and $\gep>0$.
A unit vector $v\in \mathcal H$
is called $(F,\gep)$-\emph{invariant} if
$\norm{\pi(g)v-v}<\gep$ for every $g\in F$.

A topological group $G$ has {\em property (T)} if and only if there exists a pair $(Q,\gep)$ (called a {\em Kazhdan pair}),
where $Q$ is a compact subset of $G$ and $\gep>0$,
such that for every unitary representation $\rho$ having a unit $(Q,\gep)$-invariant vector there exists a non-zero
invariant vector. Equivalently, $G$ has property (T) if and only if
the trivial representation
$1_G$ is isolated in $\mathcal R \cup \{1_G\}$ for every set $\mathcal R$
of equivalence classes of unitary representations of $G$ without non-zero invariant vectors \cite[Proposition 1.2.3]{bhv:book}.

Compact groups have property (T) \cite[Proposition 1.1.5]{bhv:book}, but countable Abelian precompact groups
do not have property (T) \cite[Theorem 6.1]{FHU}. 

We refer to Fell's papers \cite{fel1,fel2}, the classical text by Dixmier \cite{dixmier} and the recent monographs by
de la Harpe and Valette \cite{harp_vale}, and Bekka, de la Harpe and Valette
\cite{bhv:book} for basic definitions and results concerning Fell topologies and property (T).

\section{Almost metrizable groups}
\label{s:prec}

If $A$ is a subset of a topological space $X$, the {\em character} $\chi(A,X)$ of $A$ in $X$ is the least cardinality
of a base of neighborhoods of $A$ in $X$. (If this definition leads to a finite value of $\chi(A,X)$, we replace it by $\o$,
the first infinite cardinal, and similarly for other cardinal invariants.)
If $A$ is a closed subset of a compact space $X$, the character $\chi(A,X)$
equals the {\em pseudocharacter} $\psi(A,X)$ -- the least cardinality of a family $\g$ of open subsets of $X$ such that 
$\cap\g=A$. In particular, if $A$ is a closed $G_\d$-subset of a compact space $X$, then $\chi(A,X)=\o$. 

If $K$ is a compact subgroup of a topological group, then $G/K$ is metrizable if and only if $\chi(K,G)=\o$ \cite[Lemma 4.3.19]{ArhTk}.
Let $G$ be an almost metrizable topological group, $\sK$ the collection of all compact subgroups $K\sbs G$ such that $\chi(K,G)=\o$.
Then for every neighborhood $O$ of the neutral element there is $K\in\sK$ such that $K\sbs O$ \cite[Proposition 4.3.11]{ArhTk}.
We now show that if $G$ is additionally $\o$-narrow, then $K$ can be chosen normal (in the algebraic sense).  
Recall that a topological group $G$ is {\em $\o$-narrow} \cite{ArhTk} if for every neighborhood
$U$ of the neutral element there exists a countable set $A\sbs G$ such that $AU=G$.


\blem
\label{p:alm}
Let $G$ be an $\o$-narrow almost metrizable group, $\sN$ the collection of all normal (= invariant under inner automorphisms) compact
subgroups $K$ of $G$ such that the quotient group $G/K$ is metrizable (equivalently, $\chi(K,G)=\o$). Then for every neighborhood $O$ of the neutral
element there exists $K\in \sN$ such that $K\sbs O$.
\elem

\bpf
Let $L\sbs O$ be a compact subgroup of $G$ such that the quotient space $G/L=\{xL:x\in G\}$ is metrizable.
It suffices to prove that $K=\cap\{gLg\obr:g\in G\}$, the largest normal subgroup of $G$ contained in $L$,
belongs to $\sN$. 

There exists a compatible metric on $G/L$ which is invariant under the action of $G$ by left translations.
To construct such a metric, consider a countable base $U_1, U_2, \dots$ of neighborhoods of $L$ in $G$.
We may assume that for each $n$ we have $U_n=U_n\obr=U_nL$ and $U_{n+1}^2\sbs U_n$. Let 
$\g_n=\{gU_n:g\in G\}$. The open cover $\g_n$ of $G$ is invariant under left $G$-translations and under right $L$-translations,
and $\g_{n+1}$ is a barycentric refinement of $\g_n$. The pseudometric on $G$ that can be constructed in a canonical
way from the sequence $(\g_n)$ of open covers (see \cite[Theorem 8.1.10]{Eng}) gives rise to a compatible $G$-invariant
metric on $G/L$. A similar construction was used in \cite[Lemma 4.3.19]{ArhTk}.

If an $\o$-narrow group transitively acts on a metric space $X$ by isometries, then $X$ is separable \cite[10.3.2]{ArhTk}.
Thus $X=G/L$ is separable. Let $Y$ be a dense countable subset of $X$. Then
$K=\{g\in G:gx=x\hbox{ for every $x\in X$}\}=\{g\in G:gx=x\hbox{ for every $x\in Y$}\}$
is a $G_\d$-subset of $L$, hence $\chi(K,L)=\o$. It follows that $\chi(K,G)\le\chi(K,L)\chi(L,G)=\o$ (\cite[Exercise 3.1.E]{Eng}).
\epf

%

\section{Main theorem}
\label{s:main}
\bthm \label{t:1}
If $G$ is a precompact almost metrizable group, then $\widehat G$ is discrete.
\ethm

\bpf
Let $\rho$ be an irreducible unitary representation of $G$.
We must prove that $[\rho]$ is isolated in $\widehat G$. It suffices to find
a discrete open subset $D\sbs\widehat G$ such that $[\rho]\in D$.

Precompact groups are $\o$-narrow, so Proposition~\ref{p:alm} applies to $G$. 
Let $\sN$, as above, be the collection of all normal compact subgroups $K\sbs G$
such that $\chi(K,G)=\o$. Then $\sN$ is closed under countable intersections, 
and it follows from Proposition~\ref{p:alm} that for every $G_\d$-subset $A$ of $G$
containing the neutral element there exists $K\in\sN$ such that $K\sbs A$. In particular,
there exists $K\in \sN$
such that $K$ lies in the kernel of $\rho$. Let $D\sbs\widehat G$ be the set of all classes
$[\sigma]\in \widehat G$ such that $K$ is contained in the kernel of $\sigma$. Then 
$[\rho]\in D$. It suffices to verify that $D$ is open and discrete.

Step 1. We verify that $D$ is open. 
Let $\mathcal R$ be the set of equivalence
classes of all finite-dimensional unitary representations (which may be reducible) of $K$
without non-zero invariant vectors. Let $\tau_n$ be the trivial $n$-dimensional 
representation $1_K\oplus\dots\oplus1_K$ ($n$ summands) of $K$, $n=1,2\dots$. 
In the notation of section~\ref{s:prelim}, $P_{\tau_n}$ does not depend on $n$ and is 
the set of non-negative constant functions on $K$. It follows that in the space 
$\mathcal S=\mathcal R \cup \{[\tau_n]: n=1,2,\dots\}$,
equipped with the Fell topology, the points $[\tau_n]$ are indistinguishable: any open
set containing one of these points contains all the others. Since $K$ has property (T), 
$[\tau_1]=[1_K]$ is not in the closure of $\mathcal R$. Therefore $\mathcal R$ is closed in $\mathcal S$
and $\mathcal S\stm \mathcal R$ is open in $\mathcal S$.

We claim that for every irreducible unitary representation $\sigma$ of $G$ the class of the restriction
$\sigma|_K$ belongs to $\mathcal S$. In other words, the claim is that $\sigma|_K$ is trivial if it admits
a non-zero invariant vector. Let $V$ be the (finite-dimensional) space of the representation $\sigma$.
For $g\in G$ and $x\in V$ we write $gx$ instead of $\sigma(g)x$. The space $V'=\{x\in V: gx=x\hbox{ for all }g\in K\}$
of all $K$-invariant vectors is $G$-invariant. Indeed, if $x\in V'$, $g\in G$ and $h\in K$, then $g\obr hgx=x$ because
$g\obr hg\in K$ and $x$ is $K$-invariant. It follows that $hgx=gx$ which proves that $gx\in V'$. Since $\sigma$ is 
irreducible, either $V'=\{0\}$ or $V'=V$. Accordingly, either $\sigma|_K$ admits no non-zero invariant vectors or
else is trivial


We have just proved that the restriction map $r:\widehat G\to \mathcal S$ is well-defined. 
Clearly $r$ is continuous, and therefore $D=r\obr(\mathcal S\stm \mathcal R)$ is open in $\widehat G$.


Step 2. We verify that $D$ is discrete.
Let $p:G\to G/K$ be the quotient map. Then $D$ is the image of the dual map 
$\hat p:\widehat{G/K}\to \widehat G$. According to \cite[Theorem 4.1]{FHU}, 
the dual space of a metrizable precompact group is discrete. Thus $\widehat{G/K}$
is discrete. Since $p$ is a perfect map, it is compact-covering, and Lemma~\ref{l:ccov} implies that
$\hat p:\widehat{G/K}\to \widehat G$ is a homeomorphic embedding. Therefore, $D=\hat p(\widehat{G/K})$
is discrete.
\epf

\end{document}